\documentclass[reqno,11pt]{amsart}

\usepackage{amsmath,amsfonts,amssymb,amsthm,epsfig,amscd,}
\usepackage[]{fontenc}
\usepackage{enumerate}
\usepackage[cmtip,all]{xy}
\usepackage{tabularx,multirow}

 \allowdisplaybreaks \numberwithin{equation}{section}

\newtheorem{theorem}{Theorem}[section]
\newtheorem{proposition}[theorem]{Proposition}
\newtheorem{conjecture}[theorem]{Conjecture}
\newtheorem{lemma}[theorem]{Lemma}
\newtheorem{corollary}[theorem]{Corollary}
\newtheorem{claim}[theorem]{Claim}

\theoremstyle{definition}

\newtheorem{remark}[theorem]{Remark}

\newcommand{\codim}{\mathrm{codim}}

\newcommand{\Sp}{\mathrm{Sp}}

\newcommand{\Hilb}{\mathcal{H}\mathrm{ilb}}

\begin{document}

\title[]{On dominant rational maps from products of\\ curves to surfaces of general type}
\author{F. Bastianelli}
\address{Dipartimento di Matematica e Applicazioni, Universit\`a degli Studi di Milano-Bicocca, via Cozzi 53, 20125 Milano - Italy}
\email{francesco.bastianelli@unimib.it}
\thanks{This work was partially supported by INdAM (GNSAGA); PRIN 2009 \emph{``Moduli, strutture geometriche e loro applicazioni''}; FAR 2012 (PV) \emph{``Variet\`a algebriche, calcolo algebrico, grafi orientati e topologici''}}
\author{G. P. Pirola}
\address{Dipartimento di Matematica, Universit\`a degli Studi di Pavia, via Ferrata 1, 27100 Pavia - Italy}
\email{gianpietro.pirola@unipv.it}

\begin{abstract}
In this paper we investigate the existence of generically finite dominant rational maps from products of curves to surfaces of general type.
We prove that the product $C\times D$ of two distinct very general curves of genus $g\geq 7$ and $g'\geq 2$ does not admit dominant rational maps---other than the identity---on surfaces of general type.
\end{abstract}

\maketitle

\section{Introduction}\label{section INTRODUCTION}

Let $X$ be a non-singular complex projective variety of general type.
Let $\mathcal{IS}(X)$ be the \emph{Iitaka-Severi set} of $X$, whose elements are equivalence classes of generically finite dominant rational maps $X\dashrightarrow Y$ to varieties of general type, where two such maps are equivalent if they differ for a birational isomorphism $Y\dashrightarrow Y'$ of the target varieties. Moreover, let us denote by $s(X)$ its cardinality.

The finiteness of $\mathcal{IS}(X)$ has been recently proved by means of the main result of \cite{Ma} and important advances in the knowledge of pluricanonical maps (cf. \cite{Ts,Ta,HMc}).
Furthermore, there are several works providing bounds on $s(X)$ (see e.g. \cite{Ka,Tn,G,GP2,GP3} and \cite{H,NP} under additional hypothesis).

\smallskip
When $X$ is a curve of genus $g$ with general moduli (i.e. for any point $[X]\in \mathcal{M}_g$ contained in a certain non-empty Zariski open subset of the moduli space) it is well known that $X$ does not dominate other curves of genus greater than one (see e.g. \cite{Pi}). Equivalently,
\begin{theorem}\label{theorem CURVES}
Let $X$ be a general curve of genus $g\geq 2$.
Then ${s(X)=1}$.
\end{theorem}
If in addition $X$ has very general moduli (so $[X]$ lies outside some countable collection of proper subvarieties of $\mathcal{M}_g$), then its Jacobian is simple and $X$ does not dominate elliptic curves.
Furthermore, the cardinality of the Iitaka-Severi set is still one for general plane curves of degree $d\geq 4$ and for general hyperplane sections of regular surfaces (cf. \cite{CvdG,S}).
It is worth noting that these facts may be proved either by a moduli count based on Hurwitz formula, or by a monodromy argument on Hodge structures.

\smallskip
In the light of the above-mentioned results on curves, it is interesting to investigate if---under some assumption of generality in a suitable moduli space---the same statement can be extended to higher dimensional varieties of general type.

Currently, few general results are known in this direction, as it is difficult to apply the techniques used for curves to higher dimensional varieties.
On one hand, it is hard to estimate properly the number of moduli of ramification divisors of dominant rational maps $X\dashrightarrow Y$.
On the other hand, the Kodaira dimension is no longer governed by Hodge structures.

However, \cite{GP1} establishes that general surfaces $X\subset\mathbb{P}^3$ of degree $5\leq d\leq 11$ satisfy $s(X)=1$, and the same is conjectured for general hypersurfaces $X\subset \mathbb{P}^n$ of general type and arbitrary degree.

\smallskip
In the following, we discuss this issue on products of curves, that is we investigate the existence of dominant rational maps $C\times D\dashrightarrow S$, where $C$ and $D$ are smooth complex projective curves and $S$ is some surface of general type.

We recall that the product ${C\times D}$ is of general type if and only if both $C$ and $D$ are of general type.
Clearly, if either $s(C)>1$ or $s(D)>1$, then $s(C\times D)>1$.
Moreover, there are several examples of surfaces of general type dominated by products of curves, as for instance Beauville surfaces (cf. \cite[p.159]{B1} and \cite{C2,BC}).
More generally, there are series of papers aiming to classify surfaces of general type dominated by products of curves under Galois rational coverings (see e.g. \cite{BCG,BP,Pe,Po} for detailed treatises and bibliographical notes on this topic). In particular, curves involved in these constructions always posses non-trivial automorphisms.

On the other hand, when both $C$ and $D$ are assumed to have very general moduli---so they do not possess non-trivial automorphisms and maps on other curves of general type---we prove the following.
\begin{theorem}\label{theorem MAIN}
Let $C$ and $D$ be two distinct very general curves of genus $g\geq 7$ and $g'\geq 2$ respectively.
Then ${s(C\times D)=1}$.
\end{theorem}
Hence the theorem above leaves out only finitely many cases to discuss in view of a whole understanding of the problem.

This result is not a direct consequence of the one-dimensional case.
In line with \cite{GP1}, we prove Theorem \ref{theorem MAIN} by using the analogues of the techniques we mentioned above for curves.
On one hand, we exploit a Hodge theoretical argument using monodromy to deduce that the possible target surfaces of dominant rational maps $C\times D\dashrightarrow S$ have null geometric genus.
On the other hand, we perform a moduli count of products of curves, surfaces of general type having $p_g(S)=0$, and families of curves lying on them.
In particular, an important role is played by the use of Mori's bend-and-break technique (cf. \cite{Mo} and \cite[Section II.5]{Ko}) on isotrivial families of curves covering $S$, unlike the approach followed in \cite{GP1}, where it is made an important use of Hurwitz formula and bounds on number of moduli of ramification divisors.
Furthermore, in order to treat the case $g=7$, we argue by degeneration to stable curves approaching the boundary of the moduli space $\mathcal{M}_{g'}$.

\smallskip
We point out that our argument does not apply to lower genera, however we believe that it could be suitably modified to include the case of $g$ and $g'$ satisfying $g+g'\geq 9$.
Moreover, it is difficult to exhibit examples of dominant rational maps $C\times D\dashrightarrow S$ on surface of general type when $s(C)=s(D)=1$ and the curves do not posses automorphisms.
Therefore it seems natural to conjecture that Theorem \ref{theorem MAIN} could be extended both to lower genera---except when $g=g'=2$ because there always exists a dominant map $C\times D\dashrightarrow S$ (cf. Remark \ref{remark GENUS 2})---and to a Zariski open subset of $\mathcal{M}_g\times \mathcal{M}_{g'}$. Namely,
\begin{conjecture}
Let $C$ and $D$ be two distinct general curves of genus $g\geq 2$ and $g'\geq 2$ respectively.
Then
\begin{displaymath}
s(C\times D)=\left\{
\begin{array}{ll}
2 & \textrm{if }\,g=g'=2\\
1 & \textrm{otherwise}\qquad .
\end{array}\right.
\end{displaymath}
\end{conjecture}

The plan of the paper is the following.
In Section \ref{section PRELIMINARIES} we shall introduce the main preliminary results.
In Section \ref{section ISOTRIVIAL FAMILIES OF CURVES AND RIGIDITY} we shall turn to isotrivial families of curves covering surfaces of general type.
Moreover, we shall prove a rigidity result for dominant rational maps from products of curves to a fixed surface of general type.
Finally, Section \ref{section PROOF} shall be devoted to prove Theorem \ref{theorem MAIN}.

\section{Preliminaries}\label{section PRELIMINARIES}

This section concerns the preliminary results necessary to prove Theorem \ref{theorem MAIN}.
In particular, we shall recall some important facts and prove some results about the moduli space of stable curves and its boundary divisors.
Then we shall follow \cite[Section 2]{GP1} to deal with Hilbert schemes and moduli of surfaces of general type.
Finally, we shall turn to products of curves and their Hodge structures.

\smallskip
\subsection{Notation}
We work throughout over the field $\mathbb{C}$ of complex numbers.
By \emph{variety} we mean a reduced algebraic variety over $\mathbb{C}$.
When we speak of a \emph{smooth} variety, we always implicitly assume it to be irreducible.

Given a variety $X$, we say that a property holds for a \emph{general} point ${x\in X}$ if it holds on a Zariski open non-empty subset of $X$.
Moreover, we say that $x\in X$ is a \emph{very general} point if there exists a countable collection of proper subvarieties of $X$ such that $x$ is not contained in the union of those subvarieties.

As is customary, for any smooth surface $X$, we denote by ${q(X)=\dim H^{1,0}(X)}$ the \emph{irregularity} and by ${p_g(X)=\dim H^{2,0}(X)}$ the \emph{geometric genus}.

\smallskip
\subsection{Moduli of curves}
Let $\mathcal{M}_g$ be the moduli space of smooth projective curves of genus $g\geq 2$, which is an irreducible variety of dimension $3g-3$.
Let $\overline{\mathcal{M}}_g$ denote the Deligne-Mumford compactification and let $\Delta=\overline{\mathcal{M}}_g-\mathcal{M}_g$ be the boundary divisor.
We recall that the irreducible components $\Delta_i$ of $\Delta$ are closures of loci of curves with one node. In particular, the general curve of $\Delta_0$ is irreducible, whereas for $i=1,\ldots,\left\lfloor g/2\right\rfloor$, the general curve of $\Delta_i$ splits into two irreducible curves of genus $i$ and $g-i$ attached at one point (see for instance \cite[Chapter XII]{ACG} and \cite[Chapter 2]{HM}).

We are interested in closed subvarieties of $\overline{\mathcal{M}}_g$ intersecting $\Delta$ at some point representing curves with only rational and elliptic components.
When $g= 2$, it is well known that the only complete subvarieties of $\mathcal{M}_2$ are points (cf. \cite{D}), hence any positive dimensional closed subvariety of $\overline{\mathcal{M}}_2$ parameterizes also curves having only rational and elliptic components.

As far as higher genera are concerned, the following holds (cf. \cite[Corollary 2.2]{F}) .
\begin{theorem}\label{theorem FABER}
For $g=3$ and $g\geq 5$, any two irreducible closed subvarieties of $\overline{\mathcal{M}}_g$ of codimension one have non-empty intersection.
\end{theorem}
Therefore codimension-one subvarieties of $\mathcal{M}_g$ satisfy the following property.
\begin{corollary}\label{corollary DELTA_1}
Let $g\geq 3$ and let $\mathcal{Z}\subset \mathcal{M}_g$ be an irreducible subvariety of codimension one.
Then the closure $\overline{\mathcal{Z}}\subset \overline{\mathcal{M}}_g$ meets $\Delta_1$.
\begin{proof}
In virtue of Theorem \ref{theorem FABER}, the only case to treat is $g=4$.
Consider the Chow group $\displaystyle A^1\left(\overline{\mathcal{M}}_4\right)$ of codimension-one cycles of $\overline{\mathcal{M}}_4$, and let $\lambda,\delta_0,\delta_1,\delta_2$ be the standard independent divisor classes generating $\displaystyle A^1\left(\overline{\mathcal{M}}_4\right)$ as a vector space.
Let $a\lambda +a_0\delta_0+a_1\delta_1+a_2\delta_2$ be the class of $\overline{\mathcal{Z}}$.
If $\Delta_1$ and $\overline{\mathcal{Z}}$ had empty intersection, the product of their classes would vanish in the Chow group $\displaystyle A^2\left(\overline{\mathcal{M}}_4\right)$ of codimension-two cycles.
In particular, we would have a relation $\displaystyle \left(a\lambda +a_0\delta_0 +a_1\delta_1+a_2\delta_2\right)\delta_1=0$ between the products of the standard classes, but this is impossible as the unique such a relation is $\displaystyle \left(10\lambda -\delta_0 -2\delta_1\right)\delta_2=0$ (see \cite[Section 2]{F}).
\end{proof}
\end{corollary}

Thus we deduce the following.
\begin{proposition}\label{proposition RATIONAL AND ELLIPTIC CURVES}
Let $g\geq 2$ and let $\mathcal{Z}\subset \mathcal{M}_g$ be an irreducible subvariety of codimension one.
Then there exists $[Z']\in\overline{\mathcal{Z}}$ such that $Z'$ has only rational and elliptic components.
\begin{proof}
The case $g=2$ has been discussed above, so we set $g\geq 3$ and we proceed by induction on $g$.
By Corollary \ref{corollary DELTA_1}, the closure $\overline{\mathcal{Z}}\subset \overline{\mathcal{M}}_g$ meets the boundary at $\Delta_1$.
Hence there exists an irreducible complete subvariety $\Sigma\subset\overline{\mathcal{Z}}\cap\Delta_1$ of dimension $3g-5$.
In particular, if $[Z]\in \Sigma$, then $Z$ consists of some curves $Z_1$ and $Z_{g-1}$ attached at one point, with $[Z_{1}]\in \overline{\mathcal{M}}_{1}$ and $[Z_{g-1}]\in \overline{\mathcal{M}}_{g-1}$.

Suppose that for general $[Z]\in \Sigma$, the curve $Z$ possesses a (unique) smooth component $Z_{g-1}$.
Therefore the image of the projection $\pi\colon \Sigma\dashrightarrow \mathcal{M}_{g-1}$ has dimension $\displaystyle\dim \pi\left(\Sigma\right)\geq 3g-7=\dim \mathcal{M}_{g-1}-1$.
By induction, there exists $\displaystyle [Z_{g-1}']\in \overline{\pi\left(\Sigma\right)}$ such that $Z_{g-1}'$ has only rational and elliptic components.
As $\Sigma$ is complete, there exists $[Z']\in \Sigma\subset \overline{\mathcal{Z}}$, with $Z'$ consisting of $Z_{g-1}'$ and some elliptic curve $Z_1'$ attached at one point.

On the other hand, assume that the general $[Z]\in \Sigma$ has no smooth components of genus $g-1$.
Thus there exists a fixed $0\leq i\leq \left\lfloor \frac{g-1}{2}\right\rfloor$ such that any $[Z]$ is obtained from some $[Z_{1}]\in \overline{\mathcal{M}}_{1}$ and $[Z_{g-1}]\in \Delta_i\subset \overline{\mathcal{M}}_{g-1}$.
Viceversa, by irreducibility, completeness and $\dim \Sigma= \dim \overline{\mathcal{M}}_{1} + \dim\Delta_i + 1$, we deduce that for any $[Z_{1}]\in \overline{\mathcal{M}}_{1}$ and $[Z_{g-1}]\in \Delta_i$, there exists $[Z]\in \Sigma$ given by $[Z_1]$ and $[Z_{g-1}]$.
In particular, $\Delta_i$---and consequently $\Sigma$---parameterizes also curves having only rational and elliptic component.
\end{proof}
\end{proposition}

\smallskip
\subsection{Families of curves}
Let $B$ be a smooth variety and let  $\displaystyle \mathcal{E}\stackrel{q}{\longrightarrow} B$ be a family of curves of genus $g$, that is a surjective proper morphism such that $E_b=q^{-1}(b)$ is a curve of genus $g$.
We recall that if $q$ is a smooth morphism, it is naturally defined a modular map ${\mu\colon B\longrightarrow \mathcal{M}_g}$ as $\mu(b)=[E_b]$ (see e.g. \cite[Chapter 5]{MF}).
More generally, if $q$ is not smooth and $\mathcal{E}=\bigcup \mathcal{E}_i$ is an irreducible decomposition such that any component dominates $B$, we can make a base change
\begin{equation*}
\xymatrix{ \mathcal{F}_i \ar[r]^\nu \ar[d]_{p_i} & \mathcal{E}_i\ar[d]^{q} \\  W_i \ar[r] & B}
\end{equation*}
such that $\displaystyle \mathcal{F}_i\stackrel{p_i}{\longrightarrow} W_i$ is a smooth family of genus $g$ curves and ${\nu^{-1}(E_b)\longrightarrow E_b}$ is the normalization map (see \cite[Section 2.3]{GP1}).
Thus we still define the modular map ${\mu\colon W_i\longrightarrow \mathcal{M}_g}$ such that ${\mu(w)=[F_w]}$, and we define the \emph{modular dimension} of the family $\displaystyle\mathcal{E}\stackrel{q}{\longrightarrow} B$ as
\begin{equation*}
\displaystyle M(\mathcal{E}/B):=\max_i \dim \mu(W_i).
\end{equation*}

Given a surface $S$ of general type, let us denote by $\Hilb(S)$ its Hilbert scheme.
We assume further that $\displaystyle\mathcal{E}\stackrel{q}{\longrightarrow} B$ is a family of curves on $S$, that is $\displaystyle \mathcal{E}\subset B\times S$.
Over a Zariski open subset $U\subset B$, the morphism $q$ is flat and we can define a map ${\rho\colon U\longrightarrow \Hilb(S)}$ sending $b\in U$ to the point parameterizing the curve $E_b$ on $S$.
Thus we define the \emph{dimension of the family} $\mathcal{E}$ \emph{on} $S$ as
\begin{equation*}
D(\mathcal{E}/B):=\dim \rho(U),
\end{equation*}
and we have the following result based on Mori's bend-and-break (cf. \cite[Proposition 2.3.2]{GP1}).
\begin{theorem}\label{theorem B&B}
Let $\mathcal{E}\longrightarrow B$ be a family of curves of genus $g\geq 2$ on a surface of general type. Then
\begin{equation*}
M(\mathcal{E}/B)\leq D(\mathcal{E}/B)\leq M(\mathcal{E}/B)+1.
\end{equation*}
\end{theorem}

\begin{remark}\label{remark RATIONAL CURVES}
We recall a well-known fact often involved in the proof of Theorem \ref{theorem MAIN}, that is neither rational nor elliptic curves cover a surface of general type (see \cite[Proposition VII.2.1]{BHPV}).
In particular, $D(\mathcal{E}/B)=0$ for any family of curves on $S$ of genus $g<2$.
\end{remark}

\smallskip
\subsection{Moduli of surface of general type}
Let $\mathcal{M}_{K^2,\chi}$ be the variety  parameterizing isomorphism classes of surfaces having numerical invariants $\chi$ and $K^2$.
As in \cite{C1}, we consider the isomorphism class $[S]\in\mathcal{M}_{K^2,\chi}$ of a surface of general type $S$, and we denote by $\mathcal{M}$ the union of the components of $\mathcal{M}_{K^2,\chi}$ whose points are isomorphism classes of surfaces orientedly homeomorphic to $[S]$.
Then the \emph{number of moduli} $M(S)$ of $S$ is defined as the dimension of $\mathcal{M}$ at $[S]$, and it satisfies the following estimate (cf. \cite[Theorem 2.5.1]{GP1}).
\begin{theorem}\label{theorem M(S)}
Let $S$ be a minimal surface of general type. Then
\begin{equation*}
M(S)\leq 11\chi(\mathcal{O}_S)+K^2_S.
\end{equation*}
\end{theorem}

\smallskip
\subsection{Cohomology of products of curves}\label{subsection COHOMOLOGY}
Let $C$ and $D$ be two curves of genus $g$ and $g'$, respectively. Then the complex cohomology of the surface $C\times D$ is governed by \emph{K\"unneth formula} (cf. \cite[p.103-4]{GH}), that is
\begin{equation}\label{equation KUNNETH}
H^{p,q}(C\times D)\cong \bigoplus_{{i+h=p} \atop j+k=q} H^{i,j}(C)\otimes H^{h,k}(D).
\end{equation}

In particular, global canonical sections are such that $H^{2,0}(C\times D)\cong H^{1,0}(C)\otimes H^{1,0}(D)$.
Hence the canonical map ${\phi\colon C\times D\longrightarrow \mathbb{P}^{gg'-1}}$ factors as
\begin{equation}\label{equation CANONICAL MAP}
\xymatrix{ C\times D \ar[rr]^-{\phi_{|\omega_C|}\times\phi_{|\omega_D|}} & & \mathbb{P}^{g-1} \times \mathbb{P}^{g'-1}\ar[r]^-{\sigma} & \mathbb{P}^{gg'-1}},
\end{equation}
where $\phi_{|\omega_C|}$ and $\phi_{|\omega_D|}$ are the canonical maps of $C$ and $D$, whereas $\sigma$ is the Segre embedding (cf. \cite[p.87]{B2}).

We note further that $H^{2,0}(C\times D)$ is isomorphic to the space of holomorphic two-forms of the Hodge substructure ${H^1(C, \mathbb{C})\otimes H^1(D, \mathbb{C})}$, and we have the following.
\begin{lemma}\label{lemma IRREDUCIBLE}
Let $C$ and $D$ be two distinct very general curves of genus ${g\geq 2}$ and ${g'\geq 2}$, respectively.
Then the Hodge structure ${H^1(C, \mathbb{C})\otimes H^1(D, \mathbb{C})}$ is irreducible.
\begin{proof}
Let us consider the Jacobian variety $J(C)$ of $C$ and the cohomology group ${H^1(J(C), \mathbb{C})\cong H^1(C, \mathbb{C})}$.
As $C$ is assumed to be a very general curve of genus $g$, the monodromy action on ${H^1(C, \mathbb{C})}$ of the symplectic group $\Sp({2g},\mathbb{Z})$---of $2g\times 2g$ matrices preserving the intersection form on $J(C)$---is irreducible.
Analogously, we have that $\Sp({2g'},\mathbb{Z})$ acts irreducibly on ${H^1(D, \mathbb{C})}$.
Then the induced action of $\Sp({2g},\mathbb{Z})\times\Sp({2g'},\mathbb{Z})$ on ${H^1(C, \mathbb{C})\otimes H^1(D, \mathbb{C})}$ is irreducible as well (see e.g. \cite[Section I.2.7]{NS}).
\end{proof}
\end{lemma}

\section{Isotrivial families of curves and rigidity}\label{section ISOTRIVIAL FAMILIES OF CURVES AND RIGIDITY}

In this section we deal with isotrivial families of curves dominating a surface $S$ of general type.
We prove a rigidity result on dominant rational maps $C\times D\dashrightarrow S$ from a product of curves, when the first factor deforms in a family $\mathcal{C}\longrightarrow T$ with $C_0=C$.

\smallskip
We say that a family of curves $\mathcal{E}\stackrel{q}{\longrightarrow} B$ is \emph{isotrivial} if there exists a Zariski open subset $U\subset B$ such that the fibres $E_b=q^{-1}(b)$ with $b\in U$ are all isomorphic to a fixed smooth curve.
Moreover, such a family is called \emph{trivial} when it is birational to the product $B\times D$ endowed with the projection $B\times D\stackrel{p}{\longrightarrow} B$.

Now, let $D$ be a smooth projective curve of genus $g\geq 2$, and let $\mathcal{E}\stackrel{q}{\longrightarrow} B$ be a one-dimensional isotrivial family with general fibre $E_b\cong D$.
Under this assumption, the family is dominated by a trivial family (cf. \cite{Sr}), which can be constructed using the following explicit base change (see \cite[Section 2.4]{CHM}).
Let
\begin{equation}\label{equation B_0}
B_0:=\left\{(b,\psi)\left|b\in B\textrm{ and } \psi\colon E_b\longrightarrow D \textrm{ is an isomorphism} \right.\right\}
\end{equation}
and let $B'\subset B_0$ be any connected component dominating $B$. Then the fibred product
\begin{align*}
\mathcal{E}' & := B' \times_{B} \mathcal{E} \\
 & =\left\{(b,\psi, y)\left|b\in B,\,y\in E_b\textrm{ and } \psi\colon E_b\longrightarrow D \textrm{ is an isomorphism} \right.\right\}
\end{align*}
is isomorphic to $B'\times D$ under the map $\displaystyle{(b,\psi,y)\longmapsto \left((b, \psi), \psi(y)\right)}$.
Hence we have the rational map of families
\begin{equation*}
\xymatrix{B'\times D \ar@{-->}[r]^-{\beta} \ar[d] & \mathcal{E} \ar[d]  \\ B' \ar@{-->}[r] & B \\ }
\end{equation*}
given by $\left((b, \psi), y\right)\longmapsto \psi^{-1}(y)$.
Furthermore, the fibred surface $B'\times D$ is somehow universal among trivial families dominating $\mathcal{E}$.
Namely,
\begin{lemma}\label{lemma ISOTRIVIAL FAMILY}
Let $C$ and $D$ be smooth projective curves of genus $g,g'\geq 2$.
Let $\mathcal{E}\stackrel{q}{\longrightarrow} B$ be a one-dimensional isotrivial family with general fibre isomorphic to $D$.
For any dominant rational map of families
\begin{equation*}\label{equation ISOTRIVIAL FAMILY}
\xymatrix{C\times D \ar@{-->}[r]^-{\varphi} \ar[d]_{p} & \mathcal{E} \ar[d]^{q} & \\ C \ar@{-->}[r]^h & B & \\ }
\end{equation*}
there exist a curve $B'$ and a rational map $h'\colon C\dashrightarrow B'$ such that $\varphi$ factors through $B'\times D$ as a map of families
\begin{equation*}
\xymatrix{C\times D \ar@{-->}[drr]_{\varphi} \ar@{-->}[rr]^{h'\times Id_D} & & B'\times D \ar@{-->}[d]^{\beta} & \\
 & & \mathcal{E} & \!\!\!\!\!\!\!\!\!\!\!\!\!\!\!\!\!\!\!\!\!\!\!\!\!\!\!\!\!\!\!\!\!\!\phantom{a}.\\ }
\end{equation*}
\begin{proof}
Let $x\in C$ be a general point.
The fibre $D_x=p^{-1}(x)$ coincides with the curve $\{x\}\times D$, which is a copy of $D$.
Under this identification, the restriction $\varphi_{|\{x\}\times D}\colon D\longrightarrow E_{h(x)}$ is an isomorphism.
Therefore we can define a map
\begin{equation}\label{equation h'}
\begin{array}{cccl}
h'\colon & C & \dashrightarrow & B_0\\
& x & \displaystyle \longmapsto & \left(h(x),\left(\varphi_{|\{x\}\times D}\right)^{-1}\right)
\end{array}
\end{equation}
whose image is an irreducible component of $B_0$.
Then we define $B'\subset B_0$ to be a connected component containing the image $h'(C)$.
Finally, for general $(x,y)\in C\times D$, we deduce
\begin{equation*}
\left(\beta\circ\left(h'\times Id_D\right)\right)(x,y)=\beta\left(\left(h(x),\left(\varphi_{|\{x\}\times D}\right)^{-1}\right),y\right)=\varphi_{|\{x\}\times D}(y)=\varphi(x,y)
\end{equation*}
as claimed.
\end{proof}
\end{lemma}

A family of dominant rational maps from products of curves to a surface $S$ of general type may be viewed as a map $f\colon \mathcal{C}\times_T \mathcal{D}\dashrightarrow S$, where $\mathcal{C}\longrightarrow T$ and $\mathcal{D}\longrightarrow T$ are families of smooth curves, and the restrictions $f_t\colon C_t\times D_t\dashrightarrow S$ are dominant.

When one of the factors does not move, i.e. $\mathcal{D}=T\times D$ for some smooth curve $D$, we have that $\mathcal{C}\times_T \mathcal{D}\cong\mathcal{C}\times D$. In this setting, we combine Theorem \ref{theorem B&B} with Lemma \ref{lemma ISOTRIVIAL FAMILY}, and we prove that the maps of the family $f\colon \mathcal{C}\times D\dashrightarrow S$ factor through some fixed surface $B'\times D$. Namely,
\begin{proposition}\label{proposition ISOTRIVIAL FAMILY ON S}
Let $D$ be a smooth projective curve of genus $g'\geq 2$ and let $S$ be a surface of general type.
Let $\mathcal{C}\longrightarrow T$ be a family of smooth projective curves of genus $g\geq 2$ and let $f\colon \mathcal{C}\times D\dashrightarrow S$ be a rational map such that the restrictions $f_t\colon C_t\times D\dashrightarrow S$ are dominant.
Moreover, assume that for general $x\in C_t$, the curve $\{x\}\times D$ maps birationally onto its image $f_t(\{x\}\times D)\subset S$.\\
Then there exist a curve $B'$, and for general $t\in T$, a rational map $h'_t\colon C_t\dashrightarrow B'$ such that $f_t$ factors through $B'\times D$ as
\begin{equation}\label{equation ISOTRIVIAL FAMILY 2}
\xymatrix{C_t\times D \ar@{-->}[drr]_{f_t} \ar@{-->}[rr]^{h'_t\times Id_D} & & B'\times D \ar@{-->}[d] & \\
 & & S &\!\!\!\!\!\!\!\!\!\!\!\!\!\!\!\!\!\!\!\!\!\!\!\!\!\!\!\!\!\!\!\!\!\!\phantom{a}.\\ }
\end{equation}
\begin{proof}
Let $t\in T$ be a general point.
Our aim is to construct an isotrivial family $\mathcal{E}\stackrel{q}{\longrightarrow} B$ as in Lemma \ref{lemma ISOTRIVIAL FAMILY}, being independent of $t\in T$ and fitting in a commutative diagram as
\begin{equation}\label{equation ISOTRIVIAL FAMILY 3}
\xymatrix{C_t\times D \ar@{-->}[r]^-{\varphi_t} \ar[d]_{p_t} & \mathcal{E} \ar[d]^{q} & \\ C_t \ar@{-->}[r]^{h_t} & B & \!\!\!\!\!\!\!\!\!\!\!\!\!\!\!\!\!\!\!\!\!\!\!\!\!\!\!\!\!\!\!\!\!\!\phantom{a}.\\ }
\end{equation}

We recall that the fibres $p_t^{-1}(x)=\{x\}\times D$ form a family of copies of $D$ with total space $C_t\times D$.
As $f_t$ is dominant, the images $f_t(\{x\}\times D)$ produce a one-dimensional family of curves covering $S$, which is parameterized over
some curve $B_t\subset \Hilb (S)$.

Let $B:=\left(\bigcup_{t\in T} B_t\right)_{\mathrm{red}}\subset \Hilb(S)$ and let
$\mathcal{F}\subset B\times S$ be the restriction to $B$ of the universal family over $\Hilb(S)$.
We want to prove that $B$ is a curve.
By assumption, for general $t\in T$ and $x\in C_t$, the curve $\{x\}\times D$ is birational onto its image under $f_t$, and hence the normalization of $f_t(\{x\}\times D)$ is isomorphic to $D$.
Therefore the modular dimension of $\mathcal{F}\longrightarrow B$ is $M(\mathcal{F}/B)=0$.
Furthermore, the curve $B_t$ does not depend on $t\in T$.
Indeed, if $B_t$ were deformed on the Hilbert scheme of $S$ as we vary $t\in T$, the dimension of $\mathcal{F}\longrightarrow B$ on $S$ would be $D(\mathcal{F}/B)=\dim B\geq 2$, but this contradicts Theorem \ref{theorem B&B}.
Thus $\mathcal{F}\longrightarrow B$ is a one-dimensional family of curves on $S$, whose general fibre $F_b$ has normalization $D$.

Let $h_t\colon C_t\longrightarrow B$ be the map sending $x\in C_t$ to the point $b\in B$ such that $f_t(\{x\}\times D)=F_b$.
Then we define the dominant rational map
\begin{equation*}
\begin{array}{cccl}
\psi_t\colon & C_t\times D & \dashrightarrow & \mathcal{F}\\
& (x,y) & \displaystyle \longmapsto & \left(h_t(x),f_t(x,y)\right).
\end{array}
\end{equation*}
Normalizing $\mathcal{F}$---and possibly shrinking $B$---we obtain an isotrivial family of curves
$\widetilde{\mathcal{F}}\stackrel{q}{\longrightarrow} B$ whose general fibre $\widetilde{F}_b$ is isomorphic to $D$.
Furthermore, we can lift $\psi_t\colon C_t\times D\dashrightarrow \mathcal{F}$ to a dominant rational map of families $\varphi_t\colon C_t\times D\dashrightarrow \widetilde{\mathcal{F}}$ fitting in (\ref{equation ISOTRIVIAL FAMILY 3}).
Thus Lemma \ref{lemma ISOTRIVIAL FAMILY} assures that there exist a rational map $h_t'\colon C_t\dashrightarrow B_0$ as in $(\ref{equation h'})$ and a connected component $B'\subset B_0$ containing $h'_t(C_t)$, such that
\begin{equation*}
\xymatrix{C_t\times D \ar@{-->}[drr]_{\varphi_t} \ar@{-->}[rr]^{h'_t\times Id_D} & & B'\times D \ar@{-->}[d]^{\beta} & \\
 & & \widetilde{\mathcal{F}} & \!\!\!\!\!\!\!\!\!\!\!\!\!\!\!\!\!\!\!\!\!\!\!\!\!\!\!\!\!\!\!\!\!\!\phantom{a}.\\ }
\end{equation*}
Clearly, any $h_t'\colon C_t\dashrightarrow B_0$ maps on some irreducible component of $B_0$.
Over an open subset of $T$, we can then assume that the image $h'_t(C_t)$, and consequently $B'$ is independent of $t$.

Finally, let $\pi\colon\widetilde{\mathcal{F}}\longrightarrow S$ be the map inherited from the natural projection of $\mathcal{F}\subset B\times S$.
By construction $f_t=\pi\circ \varphi_t$ as rational maps, and hence we have the following commutative diagram
\begin{equation*}
\xymatrix{C_t\times D \ar@{-->}[drr]_{f_t} \ar@{-->}[rr]^{h'_t\times Id_D} & & B'\times D \ar@{-->}[d]^{\pi\circ \beta} & \\
 & & S & \!\!\!\!\!\!\!\!\!\!\!\!\!\!\!\!\!\!\!\!\!\!\!\!\!\!\!\!\!\!\!\!\!\!\phantom{a}. }
\end{equation*}
\end{proof}
\end{proposition}

\begin{remark}
We would like to note that both Lemma \ref{lemma ISOTRIVIAL FAMILY} and Proposition \ref{proposition ISOTRIVIAL FAMILY ON S} could also be proved following the approach of \cite[Section 2]{Sr}.
\end{remark}

\begin{remark}
Under the hypothesis of the proposition, we assume in addition that for general $t\in T$ and $y\in D$, the curve $C_t\times\{y\}$ maps birationally onto its image $f_t(C_t\times\{y\})\subset S$.
Then the general $C_t$ is birational to an irreducible component of $B'$, and hence the family $\mathcal{C}\longrightarrow T$ is isotrivial.
In particular, any map $C\times D\dashrightarrow S$ is rigid under deformations as above of the first factor.
\end{remark}

\begin{remark}
More generally, it would be interesting to have a rigidity theorem for dominant rational maps $C\times D\dashrightarrow S$, which deform in families $f\colon \mathcal{C}\times_T \mathcal{D}\dashrightarrow S$.\\
For instance, by using the techniques of this paper, it is possible to prove that if $D$ has genus $g'=2$, it cannot move and the maps $f_t\colon C_t\times D\dashrightarrow S$ fit in (\ref{equation ISOTRIVIAL FAMILY 2}).\\
In particular, if both $C$ and $D$ have genus 2, then $C\times D\dashrightarrow S$ is rigid as a map from products of smooth curves to a fixed surface of general type.
\end{remark}

\section{Proof of Theorem \ref{theorem MAIN}}\label{section PROOF}

This section is devoted to prove Theorem \ref{theorem MAIN}.
We start with a preliminary lemma providing restrictions on surfaces of general type dominated by products of very general curves. Namely,

\begin{lemma}\label{lemma p_g=0}
Let $C$ and $D$ be two distinct very general curves of genus ${g> 2}$ and ${g'\geq 2}$, respectively.
Let $S$ be a minimal surface of general type and let ${f\colon C\times D\dashrightarrow S}$ be a dominant rational map of degree ${m>1}$.
Then
\begin{itemize}
  \item[(\emph{i})] $p_g(S)=0$,
  \item[(\emph{ii})] $M(S)\leq 19$.
\end{itemize}
\begin{proof}
(\emph{i})
Let
\begin{equation*}
\xymatrix{ X \ar[d]_{h} \ar[dr]^{\widetilde{f}} & \\ C\times D\ar@{-->}[r]^f & S \\ }
\end{equation*}
be a resolution of the indeterminacy locus of $f$, and let ${f^*\colon H^2(S,\mathbb{C})\longrightarrow H^2(C\times D,\mathbb{C})}$ be the Hodge structure map  defined as the composition of the pullback map ${\widetilde{f}^*\colon H^2(S,\mathbb{C})\longrightarrow H^2(X,\mathbb{C})}$ with the Gysin map ${h_*\colon H^2(X,\mathbb{C})\longrightarrow H^{2}(C\times D,\mathbb{C})}$.

Let us consider the injective morphism ${f^*_{2,0}\colon H^{2,0}(S)\longrightarrow H^{2,0}(C\times D)}$, and let us recall that ${H^{2,0}(C\times D)\cong H^{1,0}(C)\otimes H^{1,0}(D)}$ is the holomorphic part in the Hodge decomposition of ${H^1(C, \mathbb{C})\otimes H^1(D, \mathbb{C})}$.
Then Lemma \ref{lemma IRREDUCIBLE} assures that the image of the monomorphism $f^*_{2,0}$ is either trivial or the whole $H^{2,0}(C\times D)$. Thus
either $p_g(S)=0$ or $H^{2,0}(S)\cong H^{2,0}(C\times D)$.

Aiming for a contradiction, let us assume $H^{2,0}(S)\cong H^{2,0}(C\times D)$.
Therefore the canonical map $\phi$ of ${C\times D}$ factors---as a rational map---through $f$ and we have the following diagram
\begin{equation}\label{diagram FACTORIZATION}
\xymatrix{ C\times D \ar[rr]^-{\phi} \ar@{-->}[dr]_{f} & & \mathbb{P}^{gg'-1}  \\ & S\ar@{-->}[ur] & \\ }
\end{equation}

If $g'>2$, both $C$ and $D$ are embedded by their canonical maps, and hence $\phi$ is an embedding by (\ref{equation CANONICAL MAP}).
Thus we have a contradiction as $\deg f>1$.

On the other hand, suppose that $g'=2$. Then the canonical map of $D$ factors through the hyperelliptic map. Hence $\deg \phi=2$ and the canonical image  $\phi(C\times D)$ is birational to $C\times \mathbb{P}^1$.
Thus $\deg f=2$ and $S$ is birational to $C\times \mathbb{P}^1$ as well, which is still a contradiction as $S$ is of general type.

\smallskip
(\emph{ii}) Since $S$ is a minimal surface of general type, we have $\chi(\mathcal{O}_S)=1-q(S)+p_g(S)\geq 1$, and hence $\chi(\mathcal{O}_S)=1$ by the first part of the proof.
By Theorem \ref{theorem M(S)} we have $M(S)\leq 11\chi(\mathcal{O}_S)+K^2_S$ and Miyaoka-Bogomolov inequality assures that $K^2_S \leq 9\chi(\mathcal{O}_S)$. Furthermore, if $K^2_S=9$, then $S$ is rigid by Yau's theorem.
Thus the number of moduli of $S$ satisfies $M(S)\leq 19$.
\end{proof}
\end{lemma}

\begin{remark}\label{remark GENUS 2}
Let $C$ and $D$ be distinct very general curves of genus $2$, with hyperelliptic involution $i$ and $j$, respectively.
Then their product admits a dominant map $C\times D\dashrightarrow Y$ on a surface of general type having ${q(Y)=0}$ and ${p_g(Y)=4}$, where $Y$ is the quotient of $C\times D$ under the involution ${(p,q)\longmapsto (i(p), j(q))}$.
Furthermore, $Y$ is the unique surface of general type dominated by $C\times D$ having positive geometric genus.

Indeed, if $f\colon C\times D\dashrightarrow S$ were a dominant rational map on another surface of general type with $p_g(S)>0$, we would argue as in Lemma \ref{lemma p_g=0} and $f$ would fit in (\ref{diagram FACTORIZATION}).
Since the monodromy group $M(\phi)$ of the canonical map is isomorphic to $\langle i,j\rangle\cong \mathbb{Z}/2\mathbb{Z}\times \mathbb{Z}/2\mathbb{Z}$, the surfaces with positive geometric genus fitting in the diagram are just $Y$, $\mathbb{P}^1\times D$ and $C\times \mathbb{P}^1$.
\end{remark}

Moreover, we recall that a very general curve $C$ of genus $g\geq 2$ does not dominate other curves of positive genus.
So if $f\colon C\longrightarrow E$ is a non-constant morphism to a curve $E$ with normalization $\widetilde{E}$,
then either $\widetilde{E}\cong C$---and $f$ is birational---or $\widetilde{E}\cong \mathbb{P}^1$.

We now prove our result.
Firstly, we shall use rigidity of dominant rational maps to set a moduli count excluding the case $g\geq 8$.
In order to rule out the case of genus $7$, we shall study the degeneration of trivial families of curves of genus $g'$ dominating a fixed surface $S$ of general type. Arguing by contradiction and using Proposition \ref{proposition RATIONAL AND ELLIPTIC CURVES}, we shall obtain some family of curves having only rational and elliptic components, which still covers $S$.

\begin{proof}[of Theorem \ref{theorem MAIN}]
Let $C$ and $D$ be two distinct very general curves of genus $g\geq 7$ and $g'\geq 2$ respectively.
Without loss of generality, we set ${g\geq g'}$.
By contradiction, let us assume the existence of a dominant rational map ${C\times D\dashrightarrow S}$ of degree ${m>1}$ on a surface $S$ of general type and---up to consider the minimal model of $S$---let us suppose $S$ to be smooth and minimal.

We define the locus ${\mathcal{S}\subset \mathcal{M}_g\times \mathcal{M}_{g'}}$ as
\begin{equation}\label{equation S}
\mathcal{S}:=\left\{\left([X],[Y]\right)\in\mathcal{M}_g\times \mathcal{M}_{g'} \left| \exists\, X\times Y\dashrightarrow S \textrm{ dominant}\right.\right\},
\end{equation}
and let $\mathcal{R}\subset \mathcal{S}$ be an irreducible component passing through $\displaystyle \left([C],[D]\right)$, endowed with the projection maps $\pi_1\colon \mathcal{R}\longrightarrow \mathcal{M}_{g}$ and $\pi_2\colon \mathcal{R}\longrightarrow \mathcal{M}_{g'}$.
\begin{claim}\label{claim SECOND PROJECTION}
The projection map $\pi_2\colon \mathcal{R}\longrightarrow \mathcal{M}_{g'}$ is generically finite.
\begin{proof}
We consider the fibre over the very general point $[D]\in \mathcal{M}_{g'}$,
\begin{equation*}
\pi_2^{-1}\left([D]\right)=\left\{\left([X],[D]\right)\in \mathcal{R} \left| \exists\, X\times D\dashrightarrow S \textrm{ dominant}\right.\right\}.
\end{equation*}
Aiming for a contradiction, we assume $\displaystyle \dim \pi_2^{-1}\left([D]\right)>0$.
Let ${\displaystyle \mathcal{T}:=\pi_1\left(\pi_2^{-1}\left([D]\right)\right)\subset \mathcal{M}_{g}}$ and notice that it has the same dimension of $\displaystyle \pi_2^{-1}\left([D]\right)$.

Since $[C]\in \mathcal{T}$ is a very general point of $\mathcal{M}_g$, it lies on the locus $\mathcal{M}_g^0$ of curves without automorphisms other than the identity.
Let ${T:=\mathcal{T}_{|\mathcal{M}_g^0}}$ and consider the restriction $\mathcal{C}\longrightarrow T$ of the universal family over $\mathcal{M}_g^0$.
By construction, any fibre $C_t$ is a curve admitting a dominant rational map $f_t\colon C_t\times D\dashrightarrow S$.
Up to make a base change, the family of curves $\mathcal{C}\longrightarrow T$ is endowed with a dominant rational map $f\colon \mathcal{C}\times D\dashrightarrow S$ having restrictions $f_t\colon C_t\times D\dashrightarrow S$.

Moreover, we note that for general $t\in T$ and $x\in C_t$, the curve $\{x\}\times D$ maps birationally onto its image $f_t(\{x\}\times D)$.
Indeed $\{x\}\times D$ is a copy of the very general curve $D$, therefore $f_t(\{x\}\times D)$ is either birational to $D$ or a rational curve.
In the latter case, the surface $S$ of would be covered by rational curves, but this is impossible (cf. Remark \ref{remark RATIONAL CURVES}).

Therefore the family $\mathcal{C}\longrightarrow T$ fulfils the hypothesis of Proposition \ref{proposition ISOTRIVIAL FAMILY ON S}. Thus there exist a curve $B'$ and---for general $t\in T$---a non-constant map $h'_t\colon C_t\dashrightarrow B'$ such that
\begin{equation}\label{equation ISOTRIVIAL FAMILY PROOF}
\xymatrix{C_t\times D \ar@{-->}[drr]_{f_t} \ar@{-->}[rr]^{h'_t\times Id_D} & & B'\times D \ar@{-->}[d] & \\
 & & S & \!\!\!\!\!\!\!\!\!\!\!\!\!\!\!\!\!\!\!\!\!\!\!\!\!\!\!\!\!\!\!\!\!\!\phantom{a}.\\ }
\end{equation}
In particular, there is an irreducible component of $B'$ dominated by the general fibre $C_t$, and hence by the very general curve $C$.
Therefore the closure of such a component must be a rational curve.
Then the surface $S$ is covered by a family of rational curves under the map $B'\times D\dashrightarrow S$ in (\ref{equation ISOTRIVIAL FAMILY PROOF}).
Thus we have a contradiction as $S$ is of general type.
\end{proof}
\end{claim}

In virtue of Lemma \ref{lemma p_g=0} we have that ${p_g(S)=0}$ and the modular dimension of $S$ satisfies ${M(S)\leq 19}$.
Furthermore, for a general choice of $S$ among minimal surfaces of general type having ${p_g=0}$ and being dominated by a product of very general curves, we have
\begin{equation}\label{equation DIMENSION}
\dim \mathcal{R}\geq \dim \left(\mathcal{M}_g\times \mathcal{M}_{g'}\right)-M(S)\geq (3g-3)+(3g'-3)-19.
\end{equation}
Moreover, the generically finiteness of $\pi_2\colon \mathcal{R}\longrightarrow \mathcal{M}_{g'}$ assures that $\dim \mathcal{R}\leq \dim \mathcal{M}_{g'}=3g'-3$. Thus $(3g-3)-19\leq 0$ and the assertion of Theorem \ref{theorem MAIN} is proved for any $g\geq 8$ and $g'\geq 2$.

\smallskip
Then we assume $g=7$ and $2\leq g'\leq 7$.
As a consequence of inequality (\ref{equation DIMENSION}) and Claim \ref{claim SECOND PROJECTION} we deduce that the image $\mathcal{Z}:=\pi_2(\mathcal{R})\subset \mathcal{M}_{g'}$ is a subvariety of dimension $3g'-4\leq \dim \mathcal{Z}\leq 3g'-3$.
Let $\overline{\mathcal{Z}}\subset \overline{\mathcal{M}}_{g'}$ be its closure.
Therefore Proposition \ref{proposition RATIONAL AND ELLIPTIC CURVES} assures that there exists $[Z']\in\overline{\mathcal{Z}}$ such that all the irreducible components of $Z'$ are rational and elliptic.

\begin{claim}\label{claim BOUNDARY}
For any $[Z]\in \overline{\mathcal{Z}}-\mathcal{Z}$, there exist a family $\mathcal{X}\longrightarrow W$ of nodal curves of genus $g'$ and a dominant rational map $\mathcal{X}\dasharrow S$, such that the general fibre $X_w$ of the family consists of a curve birational to $Z$ and some rational components.
\begin{proof}
Let $Z$ be a curve of genus $g'$ such that ${[Z]\in \overline{\mathcal{Z}}-\mathcal{Z}}$.
Let $U$ be a disk parameterizing a family $\mathcal{Y}\longrightarrow U$ of curves such that $[Y_t]\in \mathcal{Z}$ for $t\neq 0$ and $Y_0=Z$.
By construction for any $t\neq 0$, there exists a dominant rational map ${f_t\colon X_t\times Y_t\dashrightarrow S}$, with $[X_t]\in \mathcal{M}_g$.
In particular, the images of the curves $\{x\}\times Y_t$ under the maps $f_t$ describe a family of curves covering $S$.

We focus on those curves $f_t(\{x\}\times Y_t)\subset S$ passing through a fixed general point $s\in S$.
Up to some base change, we can assume the maps $f_t$ to vary holomorphically on $U^*=U-\{0\}$.
Hence we can define a dominant rational map
\begin{equation*}\label{diagram SEMISTABLE 2}
\xymatrix{\mathcal{Y} \ar@{-->}[r]^{\xi} \ar[d] &  S & \\  U^* &  & \!\!\!\!\!\!\!\!\!\!\!\!\!\!\!\!\!\!\!\!\!\!\!\!\!\!\!\!\!\!\!\!\!\!\phantom{a},\\}
\end{equation*}
where the restrictions $\xi_{|Y_t}$ are inherited from the maps $f_t$, and the curves $\xi\left(Y_t\right)$ pass through $s\in S$.
We extend $\xi$ at $t=0$ by nodal reduction (cf. \cite[p.119]{HM})
\begin{equation*}\label{diagram SEMISTABLE 3}
\xymatrix{\mathcal{Y'} \ar[r]^{\xi'} \ar[d] & S & \\  U' &  & \!\!\!\!\!\!\!\!\!\!\!\!\!\!\!\!\!\!\!\!\!\!\!\!\!\!\!\!\!\!\!\!\!\!\phantom{a},}
\end{equation*}
so that the central fibre $Y'_0$ is a nodal curve containing a curve birational to $Z$ and the remaining components are rational.
Moreover, the image $\xi'\left(Y'_0\right)\subset S$ is a curve passing through the general point $s\in S$, and hence the assertion follows.
\end{proof}
\end{claim}

Since $[Z']\in \overline{\mathcal{Z}}-\mathcal{Z}$, there exists a family $\mathcal{X}'\longrightarrow W$ and a dominant rational map $\mathcal{X}'\dashrightarrow S$ as in Claim \ref{claim BOUNDARY}.
In particular, each component of the general fibre $X'_w$ is either rational or elliptic.
Hence the surface of general type $S$ is covered by curves of genus smaller than two.
Thus we get a contradiction and Theorem \ref{theorem MAIN} is proved.
\end{proof}

\begin{remark}
We note that when $2\leq g'\leq g\leq 6$, the subvariety $\mathcal{Z}\subset \mathcal{M}_{g'}$ has no longer $\codim \mathcal{Z}\leq 1$, and our argument does not apply.
On the other hand, we are assured that $\overline{\mathcal{Z}}$ meets the boundary $\Delta\subset \overline{\mathcal{M}}_{g'}$ when $\codim \mathcal{Z}\leq 2g'-2$ (cf. \cite{D}).
Hence one can hope to extend Theorem \ref{theorem MAIN} to lower genera by refining our techniques and using some results describing the intersection between $\overline{\mathcal{Z}}$ and some $\Delta_i\subset \overline{\mathcal{M}}_{g'}$ (see e.g. \cite{CP,Kr} for $\Delta_0$).
\end{remark}

\section*{Acknowledgements}

We would like to thank Alessandro Ghigi, Valeria Marcucci, Antonio Moro and Roberto Pignatelli for helpful suggestions.

\end{document}